\documentclass[a4paper,12pt]{amsart}
\usepackage{amssymb}
\usepackage{ifthen}
\nonstopmode \numberwithin{equation}{section}
\newtheorem{thm}{Theorem}
\newtheorem{cor}{Corollary}
\newtheorem{lem}{Lemma}

\newtheorem{conj}{Conjecture}

\theoremstyle{definition}
\newtheorem{defn}{Definition}[section]

\newtheorem{prob}[equation]{Problem}

\newenvironment{rem}{%
\bigskip
\noindent \textsl{{\sl Remark. }}}{\bigskip}
\newenvironment{rems}{%
\bigskip
\noindent \textsl{{\sl Remarks. }}}{\bigskip}

\newcounter {own}
\def\theown {\thesection       .\arabic{own}}

\newenvironment{pf}[1][]{%
 \vskip 3mm
 \noindent
 \ifthenelse{\equal{#1}{}}%
  {{\slshape Proof. }}%
  {{\slshape #1.} }%
 }%
{\qed\bigskip}

\newcounter{alphabet}
\newcounter{tmp}
\newenvironment{Thm}[1][]{\refstepcounter{alphabet}%
\bigskip%
\noindent%
{\bf Theorem \Alph{alphabet}}%
\ifthenelse{\equal{#1}{}}{}{ (#1)}%
{\bf .} \itshape}{\vskip 8pt}

\newcommand{\ID}{{\mathbb D}}
\newcommand{\IDs}{{\mathbb D^*}}
\newcommand{\IDb}{{\overline{\mathbb D}}}
\newcommand{\IDsb}{{\overline{\mathbb D^*}}}

\newcommand{\IC}{{\mathbb C}}
\newcommand{\C}{{\mathbb C}}
\newcommand{\sphere}{{\widehat{\mathbb C}}}

\newcommand{\inv}{^{-1}}

\newcommand{\dz}{{\partial}}
\newcommand{\dzb}{{\bar\partial}}

\def\be{\begin{equation}}
\def\ee{\end{equation}}

\newcommand{\bee}{\begin{enumerate}}
\newcommand{\eee}{\end{enumerate}}

\newcommand{\blem}{\begin{lem}}
\newcommand{\elem}{\end{lem}}
\newcommand{\bthm}{\begin{thm}}
\newcommand{\ethm}{\end{thm}}
\newcommand{\bcor}{\begin{cor}}
\newcommand{\ecor}{\end{cor}}
\newcommand{\beg}{\begin{examp}}
\newcommand{\eeg}{\end{examp}}
\newcommand{\begs}{\begin{examples}}
\newcommand{\eegs}{\end{examples}}
\newcommand{\bdefe}{\begin{defn}}
\newcommand{\edefe}{\end{defn}}
\newcommand{\bprob}{\begin{prob}}
\newcommand{\eprob}{\end{prob}}
\newcommand{\bei}{\begin{itemize}}
\newcommand{\eei}{\end{itemize}}

\newcommand{\bcon}{\begin{conj}}
\newcommand{\econ}{\end{conj}}
\newcommand{\bcons}{\begin{conjs}}
\newcommand{\econs}{\end{conjs}}
\newcommand{\bprop}{\begin{propo}}
\newcommand{\eprop}{\end{propo}}
\newcommand{\br}{\begin{rem}}
\newcommand{\er}{\end{rem}}
\newcommand{\brs}{\begin{rems}}
\newcommand{\ers}{\end{rems}}
\newcommand{\bo}{\begin{obser}}
\newcommand{\eo}{\end{obser}}
\newcommand{\bos}{\begin{obsers}}
\newcommand{\eos}{\end{obsers}}
\newcommand{\bpf}{\begin{pf}}
\newcommand{\epf}{\end{pf}}
\newcommand{\ba}{\begin{array}}
\newcommand{\ea}{\end{array}}
\newcommand{\beq}{\begin{eqnarray}}
\newcommand{\beqq}{\begin{eqnarray*}}
\newcommand{\eeq}{\end{eqnarray}}
\newcommand{\eeqq}{\end{eqnarray*}}

\newcounter{minutes}
\divide\time by 60
\newcounter{hours}
\multiply\time by 60 \addtocounter{minutes}{-\time}

\begin{document}
\bibliographystyle{amsplain}
\title{Area distortion under meromorphic mappings with nonzero pole having quasiconformal extension}
\begin{center}
{\tiny \texttt{FILE:~\jobname .tex,
        printed: \number\year-\number\month-\number\day,
        \thehours.\ifnum\theminutes<10{0}\fi\theminutes}
}
\end{center}

\author{Bappaditya Bhowmik${}^{~\mathbf{*}}$}
\address{Bappaditya Bhowmik, Department of Mathematics,
Indian Institute of Technology Kharagpur, Kharagpur - 721302, India.}
\email{bappaditya@maths.iitkgp.ernet.in}
\author{Goutam Satpati}
\address{Goutam Satpati, Department of Mathematics,
Indian Institute of Technology Kharagpur, Kharagpur - 721302, India.}
\email{gsatpati@iitkgp.ac.in}

\keywords{ Meromorphic, Quasiconformal, Area Distortion.\\
${}^{\mathbf{*}}$ Corresponding author}
\date{ \today ; File: BS.tex}

\begin{abstract}
Let $\Sigma_k(p)$ be the class of univalent meromorphic functions defined on
$\ID$ with 
$k$-quasiconformal extension to the extended complex plane $\sphere$, where
$0\leq k < 1$. Let $\Sigma_k^0(p)$ be the class of functions $f \in \Sigma_k(p)$ having expansion of the form $f(z)= 1/(z-p) + \sum_{n=1}^{\infty}b_n z^{n}$
on $\ID.$ In this article, we obtain sharp area distortion and weighted area distortion inequalities for functions in $\Sigma_k^0(p)$. As a consequence
of the obtained results, we present a sharp estimate for the bound of the
Hilbert transform.
\end{abstract}

\maketitle
\pagestyle{myheadings}
\markboth{Area Distortion}{Area Distortion}
\bigskip

\section{Introduction}
Let $\IC$ denote the complex plane and $\sphere$ be the extended complex plane
$\C\cup\{\infty\}.$ Throughout the discussion in this article, we shall use the following notations:
$\ID=\{z : |z|<1\}$, $\overline{\ID}=\{z :
|z|\leq 1\}$, $\ID^{*}= \{z : |z|>1\}$, $\overline{\ID^{*}}= \{z :
|z|\geq1\}$.
Let $\Sigma$ be the class of
univalent meromorphic functions defined on $\ID$ having simple pole at the origin
with residue 1 and therefore each $f\in \Sigma$  has the following expansion
\begin{equation}\label{eq000}
f(z)=z^{-1} + \sum_{n=0}^{\infty}b_n z^{n}, \quad z\in \ID.
\end{equation}
It is well-known that the
univalent functions defined in $\ID$ that admit a quasiconformal
extension to the sphere $\sphere$ play an important role in Teichm\"{u}ller space theory. It is therefore of interest to study such class of functions.
To this end, let  $\Sigma_k$ be the class of functions  in $\Sigma$ that have
$k$-quasiconformal extension $(0 \leq k <1)$ to the whole plane $\sphere$.
Here, a mapping $f:\sphere\to\sphere$ is called $k$-quasiconformal
if $f$ is a homeomorphism and has locally $L^2$-derivatives on
$\IC\setminus\{f\inv(\infty)\}$
(in the sense of distribution) satisfying $|\dzb f|\le k|\dz f|$ a.e.,
where $\dz f=\partial f/\partial z$ and $\dzb f=\partial f/\partial\bar z.$
Note that such an $f$ is also called $K$-quasiconformal more often,
where $K=(1+k)/(1-k)\ge1$.
The quantity $\mu=\dzb f/\dz f$ is called the complex dilatation of $f$. The functions in the class $\Sigma_k$ has primarily been studied by O. Lehto,
(compare \cite{Lehto}) and  later R. K\"{u}hnau and S. Krushkal continued the research in this direction. More precisely,
they obtained distortion theorems, coefficient estimates, area theorem for functions in this class.

In 1955, Bojarski considered the  area distortion problem for quasiconformal mappings ( see f.i. \cite{BJ}).
Thereafter further improvements on this problem were made by Gehring and Reich (compare \cite[Theorem 1]{GR})
in a more precise form and they conjectured that

\begin{Thm}
If $f:\ID \to \ID$ be a $k$-quasiconformal mapping with $f(0)=0$, then
$$
|f(E)|\leq M(K)|E|^{1/K},
$$
for all measurable set $E \subset \ID$, where $|\cdot|$ stands for the area,
$K=(1+k)/(1-k)\geq 1,$ and the constant $M(K)=1+O(K-1),$ as $K \to 1.$
\end{Thm}
This conjecture was proved by K. Astala (\cite[Theorem~1.1]{KA}) in 1994 using thermodynamic formalism and holomorphic motion theory.
Later, Eremenko and Hamilton in \cite[Theorem 1]{EH} gave a direct and much more simpler proof to the above problem.
They assumed $f$ to be a $k$-quasiconformal mapping of the plane which is conformal on $\IC \setminus \Delta,$ where
$\Delta$ is a compact set of transfinite diameter 1 and $f$ has the normalization
$f(z)=z+o(1)$ near $\infty.$
Here we introduce the class $\Sigma_k^0$ that consists of functions
defined on $\ID^*$, having $k$-quasiconformal extension in $\ID$ such that they have pole at the point $z= \infty$ and have
the following form
\begin{equation}\label{eq00}
f(z)=z + \sum_{n=1}^{\infty}b_n z^{-n}, \quad z\in \ID^*.
\end{equation}
In \cite[Theorem 1]{EH} , if we assume $\Delta = \IDb$, then $f\in \Sigma_k^0$. We state this result below:

\begin{Thm}\label{thB}
Let $f\in \Sigma_k^0$ having the expansion of the form \eqref{eq00},
so that $f(z)-z \to 0$ as $z \to \infty.$
\begin{itemize}
\item[(i)] If $f$ is conformal on $E \subset \ID,$ then
\begin{equation*}
|f(E)| \leq  \pi^{1-1/K} |E|^{1/K}.
\end{equation*}
\item[(ii)] If $f$ is conformal outside $E \subset \ID$, then
\begin{equation*}
|f(E)| \leq K|E|.
\end{equation*}
\item[(iii)] Hence, for any arbitrary subset $E$ of $\ID$,
\begin{equation*}
|f(E)| \leq K \pi^{1-1/K} |E|^{1/K}.
\end{equation*}
\end{itemize}
All the constants in the above inequalities are best possible.
\end{Thm}
In particular, equality holds in Theorem B(i) (see \cite[p. 344]{AN}) for the function
\begin{equation}\label{eq01}
f_r(z)=
 \begin{cases}
 r^{1/K-1}z, & |z|< r,\\
 z|z|^{1/K-1}, & r\leq |z|\leq 1,\\
 z, & |z|>1,
 \end{cases}
\end{equation}
where $0<r<1$ and $f$ is conformal on $E=\lbrace z:|z|< r \rbrace$.
Next, the inequality in Theorem B(ii) is sharp for the function $f_r\inv $ and
$E=\lbrace z:r^{1/K}\leq |z| \leq 1 \rbrace$ (compare \cite[p. 324]{AIM}). Also
the inequality in Theorem B(iii) is sharp as the inequalities
in  Theorem B(i) and  Theorem B(ii) are also so.
Further, Astala and Nesi proved the weighted area distortion inequality
(\cite[Theorem 1.6]{AN}), where they considered a non negative weight function $w$ defined on a measurable set $E\subset \ID.$ We state the result below:

\begin{Thm}\label{thC}
Suppose $f\in \Sigma_k^0$ having expansion of the form \eqref{eq00} and $E\subset \ID$ such that $f$ is conformal on $E$. Let $w(z) \geq 0$ be a
(measurable) weight function defined on $E$, then
$$
\pi^{1-K}\left( \int_{E} w(z)^{1/K} \,dm \right)^K \leq \int_{E} w(z)J_f(z)\,dm \leq \pi^{1-1/K} \left(\int_{E} w(z)^K \,dm \right)^{1/K}.
$$
The inequalities are sharp. Here, $dm=dxdy$ denotes the two dimensional Lebesgue measure on the plane with $z=x+iy.$
\end{Thm}
\noindent We note here that, when $w(z)=1$ for all $z\in E$, second inequality of the above theorem yields Theorem B(i).
Area distortion results for quasiconformal mappings have several consequences.
First of all it gives the precise degree of integrability of the partial derivatives of a
$K$-quasiconformal mapping. The precise regularity of quasiconformal mappings also controls the distortion of Hausdorff dimension of
a set under a $K$-quasiregular mapping. Area distortion inequality also provides sharp bounds of Hilbert transformation of characteristic function of a set lying in the domain of a quasiconformal mapping. See \cite[chap. 13, 14]{AIM}
for details.

Let $\Sigma^0(p)$ be the class of functions that are univalent, meromorphic
on $\ID$ having a simple pole at $z=p$ with residue $1$ with the following expansion

\begin{equation}\label{eq0}
f(z)= (z-p)^{-1} + \sum_{n=1}^{\infty}b_n z^{n}, \quad z\in \ID.
\end{equation}

We emphasise here that merely considering the pole of a meromorphic function at a nonzero point not only change the normalization but provide us
with the Taylor expansion of the same function inside the disc $\{z:|z|<p\}$ along with
its other Laurent expansions. In this article we consider the class $\Sigma_k^{0}(p)$
which consists of functions in $\Sigma^{0}(p)$ that have
$k$-quasiconformal extension to the whole plane $\sphere$.
Alternatively, each function in the class $\Sigma_k^0(p)$ has
the expansion of the following form
\begin{equation}\label{eq0a}
f(z)= z(1-pz)^{-1} + \sum_{n=1}^{\infty}b_n z^{-n}, \quad z\in \ID^*.
\end{equation}
This function class $\Sigma_k^{0}(p)$,
defined above has been introduced recently in \cite{BSS}. The area theorem, coefficient
estimates and distortion inequalities for this class have also been studied recently
(compare \cite{BSS}, \cite{BS}).

In this article, we prove an area distortion inequality for functions in the
class $\Sigma_k^0(p).$ This is discussed in Theorem \ref{th1} in the next section.
Further, we obtain weighted area distortion inequality for theses functions.
This is the content of the Theorem \ref{th2} in the next section.
We point out here that Theorem \ref{th1} and Theorem \ref{th2} coincide with Theorem B
and Theorem C respectively, for $p=0$, i.e. when $f\in \Sigma_k^0.$ Finally as an application of Theorem \ref{th1}, we present
a sharp estimate for the Hilbert transform of the characteristic function
$\chi_{E},$ where $E\subset \ID.$

\section{Main Results}
\noindent We start the Section with area distortion inequality for functions in the class $\Sigma_k^0(p)$.
\bthm\label{th1}
Let $f\in \Sigma_k^0(p)$ has the expansion of the form \eqref{eq0}.

\begin{itemize}
\item[(i)] If $f$ is conformal on $E \subset \ID^*,$ then
\begin{equation}\label{eq1}
|f(E)| \leq  \left[\pi(1-p^2)^{-2}\right]^{1-1/K} |f_0(E)|^{1/K}.
\end{equation}
\item[(ii)] If $f$ is conformal outside a compact set $E \subset \ID^*$, then
\begin{equation}\label{eq2}
|f(E)| \leq K|f_0(E)|.
\end{equation}
\item[(iii)] Hence, for any arbitrary subset $E$ of $\ID^*$,
\begin{equation*}
|f(E)| \leq K \left[\pi(1-p^2)^{-2}\right]^{1-1/K} |f_0(E)|^{1/K}.
\end{equation*}
\end{itemize}
Here $K=(1+k)/(1-k)$ and $f_0(z)= 1/(z-p),~ z\in \IC.$
The constants appearing in the theorem are best possible.
\ethm
\bpf (i)
Let us define $g(z):=f(1/z)$, so that $g\in \Sigma_k^0(p)$ with the
expansion of the form \eqref{eq0a} in $\ID^*$. As $g$ is obtained by composing
a M\"{o}bius transformation with a $k$-quasiconformal map $f$ in $\sphere$,
it is also $k$-quasiconformal in $\sphere$.
Here, since $f$ is conformal in $\ID$, therefore $g$ is also conformal in $\ID^*$ and hence the dilatation of $g$
has support in $\IDb$ and it has the same modulus as that of $f$.
Since $f$ is conformal on $E\subset \ID^*$, so $g$ is again
conformal on $\tilde{g}(E)= E^{\prime} \subset \ID$, where $\tilde{g}(z)=1/z$.
As a result, the dilatation $\mu$ of $g$ satisfies $|\mu(z)|\leq k$
for all $z\in \IDb \setminus E^{\prime}$ and vanishes on $E^{\prime}$.
Now we consider the dilatation
\be\label{eq3a}
\mu_{\lambda}(z)= \frac{\lambda \mu(z)}{k},~ \lambda\in \ID.
\ee
Therefore by Measurable Riemann Mapping theorem (see \cite[p.168]{AIM}), there
exists a unique quasiconformal mapping $g(z,\lambda)=g_\lambda(z)$ (for each
$\lambda$), whose dilatation is $\mu_\lambda(z)$. Now $g_\lambda \in \Sigma_{|\lambda|}^0(p)$ as $g \in \Sigma_k^0(p)$ and also $g_\lambda$ satisfies the normalization, $g_\lambda(z)= z/(1-pz)+ o(1)$ as $z\to \infty $. Hence
$g_\lambda|_{\ID^*}\in \Sigma^0(p)$, so by Chichra's area theorem
(see \cite{CHI}), we have
\begin{equation*}
|g_\lambda(\ID)| = \pi(1-p^2)^{-2}-\pi\sum_{n=1}^{\infty} n|b_n|^2
 \leq  \pi(1-p^2)^{-2}.
\end{equation*}
Thus
$$
\int_{\ID} J_\lambda(z) \,dm \leq \pi(1-p^2)^{-2}, \quad{(z=x+iy),}
$$
where $J_\lambda$ denotes the Jacobian of the map $g_\lambda.$ As
$ E^{\prime} \subset \ID$, it follows that
\begin{equation}\label{eq4}
\int_{E^{\prime}} (1-p^2)^2 \pi^{-1} J_\lambda(z) \,dm \leq 1.
\end{equation}

Now by holomorphic dependence of the solution to the Beltrami equation,
on parameter (see f.i. \cite[II, Theorem 3.1]{LV1}), the function
$\lambda \to g(z,\lambda)$ is holomorphic in the variable $\lambda \in \ID ,$
for each fixed $z\in \ID.$ This dependency also happens for the function
$\dz g(z,\lambda)$ where $g(z,\lambda)$ is analytic in $z.$ As $g(z)$ is
conformal in $E^{\prime},$ so is $g(z,\lambda)$, hence we can say that the function $\lambda \to \dz g(z,\lambda)$ is holomorphic in $\lambda\in \ID$,
for each fixed $z \in E^{\prime}.$
Since $g(z,\lambda)$ is $|\lambda|$-quasiconformal with dilatation $\mu_\lambda(z)$ in the variable $z\in \ID$, for each fixed $\lambda$,
we can write
$$
J_\lambda(z)= |\dz g(z,\lambda)|^2 - |\dzb g(z,\lambda)|^2
= |\dz g(z,\lambda)|^2(1-|\mu_\lambda(z)|^2).
$$
Thus for $z\in E^{\prime}$ we have $J_\lambda(z)=|\dz g(z,\lambda)|^2.$
As $g(z,\lambda)$ is quasiconformal in $\ID$, the Jacobian $J_\lambda(z)$
never vanishes in $\ID$ and in particular in $E^{\prime}.$ Hence, the
function $\dz g(z,\lambda)$ is a non vanishing analytic function on
$E^{\prime}\times \ID$ and so is the function
$(1-p^2)^2 \pi^{-1}\dz g(z,\lambda)^2.$
Now if we define
\begin{equation*}
a(z,\lambda)=(1-p^2)^2 \pi^{-1}|\dz g(z,\lambda)|^2,
\end{equation*}
then $\log a(z,\lambda)$ is harmonic in $\lambda\in\ID$, for $z\in E^{\prime}$. Thus from \eqref{eq4} we see that the function $a(z,\lambda)$ satisfies
the conditions of the continuous version of Lemma 1 in \cite{EH}, consequently we have
\begin{align*}
(1-p^2)^2 \pi^{-1} \int_{E^{\prime}} |\dz g(z,\lambda)|^2\,dm ~&
\leq  \left[ (1-p^2)^2 \pi^{-1} \int_{E^{\prime}} |\dz g(z,0)|^2 \,dm \right]^
{\frac{1-|\lambda|}{1+|\lambda|}}\nonumber\\
&= \left[ (1-p^2)^2 \pi^{-1} \int_{E^{\prime}} J_0(z) \,dm \right]^
{\frac{1-|\lambda|}{1+|\lambda|}}\nonumber\\
&= \left[ (1-p^2)^2 \pi^{-1} |g_0(E^{\prime})|\right]^{\frac{1-|\lambda|}{1+|\lambda|}}.
\end{align*}
Using the fact that for $z\in E^{\prime},~J_\lambda(z)=|\dz g(z,\lambda)|^2,$
we get from the above inequality
\begin{equation*}
(1-p^2)^2 \pi^{-1}|g_\lambda(E^{\prime})|\leq \left[ (1-p^2)^2 \pi^{-1} |g_0(E^{\prime})|\right]^{\frac{1-|\lambda|}{1+|\lambda|}}.
\end{equation*}

\noindent Now for $\lambda=k$, we have $g_\lambda=g$, which yields after simplification
\be\label{eq5a}
|g(E^{\prime})| \leq  \left[\pi(1-p^2)^{-2}\right]^{1-1/K} |g_0(E^{\prime})|^{1/K}.
\ee

Now since $f(z)=g(1/z),$ we get inequality \eqref{eq1}, where
$E\subset \ID^*$ and $g_0$ is replaced by $f_0$. We now find explicitly the function $g(z,0)=g_0(z)$. For $\lambda =0,$
the function $g_0$ is conformal on the whole sphere $\sphere$ onto itself as well as it satisfies the normalization of the class $\Sigma^0(p)$ on $\ID^*$, viz.
\begin{itemize}
\item[(i)] $g_0(z)-z/(1-pz) \to 0$ as $z\to \infty$,
\item[(ii)] $g_0(1/p)= \infty$,
\item[(iii)]$(1-pz)^2g_0^\prime(z)\big|_{z=1/p}=1$.
\end{itemize}
It is now easy to see that $g_0(z)=z/(1-pz)$ for all $z\in \IC$, is the only choice and hence $f_0(z)=g_0(1/z)=1/(z-p)$ for all $z\in \IC$, which proves the theorem.

Now we consider the equality case. We observe that equality holds in (\ref{eq1}) if it does hold in \eqref{eq5a} and to establish this, we consider
the following function:
\begin{equation}\label{eq6}
g(z)=
 \begin{cases}
\frac{r^{1/K-1}}{1-p^2}  \left(\frac{z-p}{1-pz}\right) + \frac{p}{1-p^2}, & z\in B(r),\\
\frac{1}{1-p^2}  \left(\frac{z-p}{1-pz}\right) \left|\frac{z-p}{1-pz}\right|^{1/K-1} +\frac{p}{1-p^2}, & z\in \IDb \setminus B(r),\\
 \frac{z}{z-pz}, & z\in \IDs,
 \end{cases}
\end{equation}
where $0<r<1$ and $B(r)~ (\subsetneq \ID)$ is the disk given by
$$
B(r)= \left\lbrace z : \left|z-\frac{p(1-r^2)}{1-p^2r^2}\right|<\frac{r(1-p^2)}{1-p^2r^2}\right\rbrace .
$$

It is easy to verify that $g$ is a member of $\Sigma_k^0(p)$ and that $g$ is
conformal on the set $E^{\prime}=B(r)\subset \ID$. To establish the equality case, we again
observe that the M\"{o}bius transformations $(z-p)/(1-pz)$ and $g_0(z)=z/(1-pz)$
maps the above disk $B(r)$ onto the disks
$\lbrace w : |w|<r \rbrace$ and
$\lbrace w :|w - p(1-p^2)^{-1}| < r(1-p^2)^{-1}\rbrace$ respectively.
Hence the right hand side of (\ref{eq5a}) becomes $\pi r^{2/K}(1-p^2)^{-2}$.
Again $g$ in \eqref{eq6} maps the disk $B(r)$ onto the disk
$\lbrace w :|w - p(1-p^2)^{-1}| < r^{1/K}(1-p^2)^{-1}\rbrace$,
which yields $|g(B(r))|= \pi r^{2/K}(1-p^2)^{-2}$.
Hence equality holds in \eqref{eq5a} for the above $g$ and $E^{\prime}=B(r).$
Now as $f(z)=g(1/z)$, we obtain the following extremal function for the inequality
\eqref{eq1}:
\begin{equation*}
f(z)=
\begin{cases}
\frac{r^{1/K-1}}{1-p^2}  \left(\frac{1-pz}{z-p}\right) + \frac{p}{1-p^2}, & z\in
\tilde{B}(r),\\
\frac{1}{1-p^2}  \left(\frac{1-pz}{z-p}\right) \left|\frac{1-pz}{z-p}\right|^{1/K-1} +\frac{p}{1-p^2}, & z\in \IDsb \setminus \tilde{B}(r),\\
\frac{1}{z-p}, & z\in \ID,
\end{cases}
\end{equation*}
where we assume $0\leq p<r<1$. Here $\tilde{B}(r) (\subsetneq \IDs)$ is the image of the
disk $B(r)$ under the map $\tilde{g}(z)=1/z$, given by
$$
\tilde{B}(r)= \left\lbrace z \in \IC : \left|z+\frac{p(1-r^2)}{r^2-p^2}\right| >
\frac{r(1-p^2)}{r^2-p^2}\right\rbrace .
$$
Hence equality holds in \eqref{eq1} for the above $f$ and $E=\tilde{B}(r)$.

\vspace{.2cm}
(ii) As before we start the proof of this part with the transformation $g(z)=f(1/z)$.
By the given condition, $g$ is conformal outside a compact set $\tilde{g}(E)= E^{\prime} \subset \ID$, where $\tilde{g}(z)=1/z$. Thus dilatation $\mu$ of $g$ vanishes outside
the compact set $E^{\prime}$. As $g\in \Sigma_k^0(p)$ of the form \eqref{eq0a} in $\IDs$,
hence by equation (1.7) of \cite[p.3]{BS}, we get
$$
g(z)=z/(1-pz)+T[\overline{\partial}g](z).
$$
Taking partial derivative of both sides w.r.t. $z$ and using
$\partial T[\omega]=H[\omega]$, we have
\begin{equation}\label{eq7}
\partial g(z) = 1/(1-pz)^2+ H [\overline{\partial}g](z),
\end{equation}
where $`T$' and $`H$' denote two dimensional `Cauchy' and `Hilbert' transform respectively (see f.i. \cite[I \S 4.3]{LV1}).
Since $\overline{\partial}g=\mu \partial g$, the above equation takes the following form
\begin{equation}\label{eq8}
\overline{\partial}g(z) = \mu/(1-pz)^2 + \mu H [\overline{\partial}g](z).
\end{equation}
It is also known that
\be\label{eq9}
w = \overline{\partial}g
= \mu(1-pz)^{-2} + \mu H\left[\mu(1-pz)^{-2}\right] + \mu H\left[\mu H\left[\mu(1-pz)^{-2}\right]\right]+ \cdots
\ee
satisfies equation \eqref{eq8} (see \cite[p.5]{BS}). By our assumption,
$w=\dzb g$ vanishes outside $E^{\prime}$. Hence using \eqref{eq7}
and the fact that the Hilbert transform is a linear isometry on $L^2(\IC)$, we get
\beq
|g(E^{\prime})|&=& \int_{E^{\prime}}J_g(z)\,dm \nonumber\\
&=&\int_{E^{\prime}}(|\partial g|^2 - |\overline{\partial}g|^2)\,dm \nonumber\\
&=& \int_{E^{\prime}}\left(\left|(1-pz)^{-2} + H[w]\right|^2-|w|^2\right)\,dm \nonumber\\
&=& \int_{E^{\prime}} \left( |1-pz|^{-4} + 2 \mathrm{Re} \left((1-p\overline{z})^{-2}
H[w]\right)\right) dm + \int_{E^{\prime}}\left(|H[w]|^2-|w|^2\right) dm \label{eq9a} \\
&\leq & \int_{E^{\prime}}|1-pz|^{-4} \,dm + 2\int_{E^{\prime}}
\left|(1-pz)^{-2} H[w] \right| dm + \int_{\IC}\left(|H[w]|^2-|w|^2\right)
dm \nonumber\\
&=& |g_0(E^{\prime})| + 2\int_{E^{\prime}}\left| (1-pz)^{-2} H[w] \right| \,dm, \label{eq10}
\eeq
where $g_0(z)=z/(1-pz)$, as mentioned earlier. Now using the fact that Hilbert transformation is linear, we get from the identity \eqref{eq9} that
$$
(1-pz)^{-2}H[w] = (1-pz)^{-2} H \left[\mu(1-pz)^{-2}\right] + (1-pz)^{-2} H\left[\mu H\left[\mu(1-pz)^{-2}\right]\right]+\cdots .
$$
This gives
\beq\label{eq11}
\int_{E^{\prime}}\left| (1-pz)^{-2}H[w] \right| \,dm ~\leq ~  \int_{E^{\prime}} |1-pz|^{-2}\left|H\left[\mu(1-pz)^{-2}\right]\right| \,dm \\
\hspace{4cm} + \int_{E^{\prime}} |1-pz|^{-2}\left|H\left[\mu H\left[\mu(1-pz)^{-2}\right]\right]\right| \,dm +\nonumber \cdots.
\eeq
We now apply `Cauchy-Schwartz' inequality and the isometry property of Hilbert transformation to the $n$-th term of the right hand side of \eqref{eq11}
to get an upper bound for this term. We show below the computational details:
\beqq
&&\int_{E^{\prime}} |1-pz|^{-2}  \Big| \underbrace{H\,\big[\,\mu H \cdots \mu H}
_\text{$n$ terms} \,[\,\mu(1-pz)^{-2}]\,\big]\Big| \,dm\\
&\leq & \left(\int_{E^{\prime}} |1-pz|^{-4} \,dm \right)^{1/2}
\bigg(\int_{E^{\prime}} \Big|\underbrace{H\,\big[\,\mu H \cdots \mu H}
_\text{$n$ terms} \,[\,\mu(1-pz)^{-2}] \,\big]\Big|^2 \,dm \bigg)^{1/2}\\
&\leq & |g_0(E^{\prime})|^{1/2} \,\bigg(\int_{\IC} \Big|
\underbrace{H\,\big[\,\mu H \cdots \mu H}
_\text{$n$ terms} \,[\,\mu(1-pz)^{-2} ]\big]\Big|^2 \,dm \bigg)^{1/2}
\eeqq
\beqq
&= & |g_0(E^{\prime})|^{1/2} \bigg(\int_{\IC}
\Big|\underbrace{\,\mu H \, \big[\, \mu H  \cdots  \mu H}_\text{$(n-1)$ terms}\,
[\,\mu(1-pz)^{-2}] \,\big]\Big|^2 \,dm \bigg)^{1/2}\\
&\leq & \|\mu\|_{\infty} |g_0(E^{\prime})|^{1/2} \,\bigg(\int_{E^{\prime}}
\Big|\underbrace{\, H\,\big[\, \mu H \cdots  \mu H}_\text{$(n-1)$ terms}\,
[\,\mu(1-pz)^{-2} ]\,\big]\Big|^2 \,dm \bigg)^{1/2}\\
& \vdots & \\
&\leq & \|\mu\|_{\infty}^n |g_0(E^{\prime}|)^{1/2} ~\left(\int_{E^{\prime}} |1-pz|^{-4}
\,dm\right)^{1/2}\\
&=& k^n |g_0(E^{\prime})|,
\eeqq
where $\|\mu\|_{\infty} =k<1.$ Using this estimate, we get from \eqref{eq11} that
\beqq
\int_{E^{\prime}}\left| (1-pz)^{-2}H[w] \right| dm &\leq &
\sum_{n=1}^{\infty} |g_0(E^{\prime})|k^n\\
&=& k(1-k)\inv |g_0(E^{\prime})|.
\eeqq
Plugging the above estimate in \eqref{eq10}, we finally obtain
\be\label{eq11a}
|g(E^{\prime})| \leq \left(\frac{1+k}{1-k}\right) |g_0(E^{\prime})|=K|g_0(E^{\prime})|.
\ee

Now applying $f(z)=g(1/z),$ we get inequality \eqref{eq2},
where $E\subset \IDs$ and $f_0(z)=1/(z-p), z \in  \IC$.
Next we show that the constant `$K$' in Theorem \ref{th1}(ii) is best possible. This
can be verified if we can show that the constant `$K$' in \eqref{eq11a} is best possible.
We consider the following example:
\begin{equation}\label{eq12}
h(z)=
 \begin{cases}
\frac{r^{1-1/K}}{1-p^2} \left(\frac{z-p}{1-pz}\right) + \frac{p}{1-p^2}, & z\in B_0(r),\\
\frac{1}{1-p^2}  \left(\frac{z-p}{1-pz}\right)\left|\frac{z-p}{1-pz}\right|^{K-1}+\frac{p}{1-p^2}, & z\in \IDb \setminus B_0(r),\\
 \frac{z}{z-pz}, & z\in \IDs,
 \end{cases}
\end{equation}
where $B_0(r)~(\subsetneq \ID)$ is the disk given by
$$
B_0(r)= \left\lbrace z : \left|z-\frac{p(1-r^{2/K})}{1-p^2r^{2/K}}\right|<\frac{r^{1/K}
(1-p^2)}{1-p^2r^{2/K}}\right\rbrace, \quad\text{for}\,\, 0<r<1.
$$


 As similar to example \eqref{eq6}, the functions $z/(1-pz)(=g_0(z))$  and  $(z-p)/(1-pz)$ maps the disk $B_0(r)$ onto the disks
$\lbrace w :|w - p(1-p^2)^{-1}| < r^{1/K}(1-p^2)^{-1}\rbrace$ and
$\lbrace w : |w|<r^{1/K}\rbrace$ respectively. This in turn implies
$|g_0(B_0(r))| = \pi r^{2/K}(1-p^2)^{-2}$
and that the function $h$ in \eqref{eq12} itself maps the disk $B_0(r)$
onto the disk
$\lbrace w :|w - p(1-p^2)^{-1}| < r(1-p^2)^{-1}\rbrace$.
To verify the assertion we set `$E^{\prime}$' in this case, as $E^{\prime}=\IDb \setminus B_0(r)$.
Then $h$ is conformal on outside of the compact set $E^{\prime}$ and
$$
|g_0(E^{\prime})|=|g_0(\IDb)| - |g_0(B_0(r))| = \pi (1-p^2)^{-2}(1-r^{2/K}).
$$
On the other hand,
\beqq
|h(E^{\prime})| &=& |h(\IDb)| - |h(B_0(r))| \\
&=& \pi (1-p^2)^{-2}(1-r^2)\\
&=& \pi (1-p^2)^{-2} -  \pi (1-p^2)^{-2}\left[1-(1-r^{2/K})\right]^K \\
&=& \pi (1-p^2)^{-2} \left[ K(1-r^{2/K}) - (K/2)(K-1)(1-r^{2/K})^2 + \cdots \right]\\
&=& K|g_0(E^{\prime})| + O\left(|g_0(E^{\prime})|^2\right), \quad\text{as} ~~~~
|g_0(E^{\prime})|\to 0.
\eeqq
Hence the constant `$K$' can not be improved as equality holds in \eqref{eq11a}
for $|g_0(E^{\prime})|$ small enough. Composing $h$ with the inverse mapping
$\tilde{g}(z)=1/z$ and taking inversion of the disk $B_0(r)$ (for $p<r$),
extremality of \eqref{eq2} follows easily, as similar to Theorem \ref{th1}(i).

\vspace{.2cm}
\noindent (iii) To prove the last part of the theorem, we consider the following change of
variable $g(z)=f(1/z)$.
Hence $g \in \Sigma_k^0(p)$ such that it is conformal on $\IDs$ and
$k$-quasiconformal on $\IDb$.
We write $g=g_1\circ g_2,$ where $g_2$ is conformal on $E \subset \ID$,
$k$-quasiconformal on $\IDb \setminus E$ and $g_2 \in \Sigma_k^0(p).$
We assume that the function $g_1$ is $k$-quasiconformal on $g_2(E)$ and hence on
$g_2(\overline{E})$ (as a set of area zero is removable for quasiconformality),
so that $g_1$ is conformal outside the compact set
$g_2(\overline{E})$ and satisfies the conditions of Theorem B(ii).
Applying Theorem \ref{th1}(i)
to $g_2$ and Theorem B(ii) to $g_1$, we get
$$
|g(E)|=|g_1(g_2(E))|\leq K|g_2(E)| \leq K \left[\pi(1-p^2)^{-2} \right]^{1-1/K}
|g_0(E)|^{1/K}.
$$
Putting $f(z)=g(1/z)$ we obtain the theorem in terms of $f$ and $g_0$ is replaced by
$f_0(z)=1/(z-p).$ As the constants in corresponding theorems for $g_1$ and $g_2$
are best possible, hence for Theorem~\ref{th1}(iii) also.
\epf

\br
For the case $p=0$, i.e. whenever $f \in \Sigma_k^0$, the
inequality \eqref{eq5a} reduces to that of Theorem B(i),
and the extremal function $g$ defined in \eqref{eq6} becomes
$f_r$, as defined in \eqref{eq01}. This coincidence
also occurs for Theorem \ref{th1}(ii), when $p=0$, as can be seen from the
inequality \eqref{eq11a} and the extremal function $h$ defined in \eqref{eq12}.
In this case $h$ reduces to $f_r \inv$ for $p=0$, which is the extremal case for
Theorem B(ii). Although, in our case $h$ is not the inverse mapping of $g$.
\er

Next we consider the weighted area distortion problem for a function in the class
$\Sigma_k^0(p)$, where we consider a nonnegative weight function $w$ defined on a subset $E$ of $\IDs$.
\bthm\label{th2}
Suppose $f\in \Sigma_k^0(p)$ with the expansion of the form \eqref{eq0}
and $E\subset \IDs$, such that $f$ is conformal on
$E$. Let $w(z) \geq 0$ be a (measurable) weight function defined on $E$, then
\beq\label{eq13}
\left[\frac{\pi}{(1-p^2)^2} \right]^{1-K} \left(\int_{E} w(z)^{1/K} J_0(z)\,dm\right)^K  \leq  \int_{E} w(z)J_f(z)\,dm \\
\leq  \left[\frac{\pi}{(1-p^2)^2}\right]^{1-1/K} \left(\int_{E} w(z)^K J_0(z)\,dm\right)^{1/K},
\nonumber
\eeq
where $J_f$ and $J_0$ denotes Jacobian of the function $f$ and
$f_0(z)=1/(z-p), z\in \IC$ respectively. The inequalities are sharp.
\ethm
\bpf
The case $w(z)=0$ for all $z$ is trivial. So we assume $w(z)>0$ for all $z \in E.$
To establish the theorem we follow the lines of the proof of \cite[Theorem 1.6]{AN}.
For the sake of completeness, we provide computational details.
Let $g(z)=f(1/z)$ having expansion of the form \eqref{eq0a} in $\IDs$. Next we consider the weight function $w_0(z)=w(1/z)$ defined on $\tilde{g}(E)= E^{\prime} \subset \ID$,
where $\tilde{g}(z)=1/z$. Therefore $g$ is conformal on $E^{\prime}$ and
$k$-quasiconformal on $\IDb \setminus E^{\prime}.$
As similar to \eqref{eq3a}, we consider the function $g_{\lambda}(z)$ with the
dilatation $\lambda k^{-1}\mu(z)$ for $\lambda\in \ID$. Again $g_{\lambda}(z)$ is
conformal on $E^{\prime}$ (since $g$ is so) and
\be\label{eq14}
g_{\lambda}^{\prime}(z) \neq 0,\quad\text{for all} \,\, z\in E^{\prime}\quad\text{and}
\,\, \lambda\in \ID.
\ee
Using the concavity of logarithm and `Jensen's Inequality', we get for any function
$a(z)>0$  defined in $E^{\prime}$, that
\be\label{eq15}
\log \left( \int_{E^{\prime}} a(z)\,dm \right) = \sup\limits_{q(z)}\left[ \int_{E^{\prime}} q(z) \log\left(\frac{a(z)}{q(z)}\right)\,dm\right],
\ee
where the supremum is taken over all functions $q(z)$ defined on $E^{\prime}$, such that\\
(i) $ 0<q(z)<1$, a.e. $z\in E^{\prime}$ and (ii) $\int_{E^{\prime}} q(z)\,dm=1.$
In our case, we take
$$
a(z)= (1-p^2)^2\pi^{-1}w_0(z)J_{\lambda}(z)=
(1-p^2)^2\pi^{-1}w_0(z)|g_{\lambda}^{\prime}(z)|^2,~~~~z\in E^{\prime},
$$
since for $z\in E^{\prime},J_{\lambda}(z)=|\dz g_{\lambda}(z)|^2=
|g_{\lambda}^{\prime}(z)|^2.$
Hence using \eqref{eq15}, we get
\beq\label{eq16}
&& \log \left(\int_{E^{\prime}} (1-p^2)^2\pi^{-1} w_0(z)|g_{\lambda}^{\prime}(z)|^2\,dm\right)\nonumber\\
&=& \sup\limits_{q(z)}\left[ \int_{E^{\prime}} q(z)\log\left(\frac{(1-p^2)^2 \pi^{-1} w_0(z)|g_{\lambda}^{\prime}(z)|^2}{q(z)} \right) \,dm \right] \nonumber\\
&=& \sup\limits_{q(z)}\left[ \int_{E^{\prime}}q(z)\log (w_0(z))\,dm + h_p(\lambda) \right],
\eeq
where
$$
h_p(\lambda)= \int_{E^{\prime}} q(z)\log\left(\frac{(1-p^2)^2 \pi^{-1}
|g_{\lambda}^{\prime}(z)|^2}{q(z)} \right) \,dm
$$
is harmonic in $\lambda \in \ID$, by \eqref{eq14}, for each $z\in E^{\prime}$. Using \eqref{eq15} and \eqref{eq4} successively, we get
$$
h_p({\lambda}) \leq \log \left( \int_{E^{\prime}} (1-p^2)^2 \pi^{-1}|g_{\lambda}^{\prime}(z)|^2 \, dm \right)\leq 0.
$$
So for each $z\in E^{\prime},~ h_p(\lambda)$ is harmonic and nonpositive in $\ID$.
Hence by using `Harnack's Inequality' and the fact that $g_0(z)=z/(1-pz)$ (as claimed in the proof of Theorem \ref{th1}(i)), we have
\beqq
h_p(\lambda) &\leq & (1-|\lambda|)(1+|\lambda|)^{-1}h_p(0)\\
&=& (1-|\lambda|)(1+|\lambda|)^{-1} \int_{E^{\prime}} q(z)\log\left(\frac{(1-p^2)^2 \pi^{-1} |g_{0}^{\prime}(z)|^2}{q(z)} \right) \,dm.
\eeqq
For $\lambda=k$, we have $g_{\lambda}=g$ and $(1+k)/(1-k)=K$. Thus using above inequality (for $\lambda =k$) in \eqref{eq16}, and also using \eqref{eq15}
once more, we get
\beqq
&& \log \left(\int_{E^{\prime}} (1-p^2)^2\pi^{-1} w_0(z)J_{g}(z)\,dm\right)\\
&\leq &  \sup\limits_{q(z)}\left[ \int_{E^{\prime}}q(z)\log w_0(z)\,dm + \frac{1}{K} \int_{E^{\prime}} q(z)\log\left(\frac{(1-p^2)^2 \pi^{-1} J_{g_{0}}(z)}{q(z)} \right) \,dm  \right]\\
&=& \frac{1}{K}\sup\limits_{q(z)}\left[\int_{E^{\prime}} q(z)\log\left(\frac{(1-p^2)^2 \pi^{-1}w_0(z)^K J_{g_{0}}(z)}{q(z)} \right) \,dm   \right]\\
&=& \log \left(\int_{E^{\prime}} (1-p^2)^2 \pi^{-1}w_0(z)^K J_{g_{0}}(z)\,dm
\right)^{1/K}.
\eeqq
Taking exponentiation and doing a rearrangement, we obtain
\be\label{eq16a}
\int_{E^{\prime}} w_0(z)J_g(z)\,dm
\leq  \left[\frac{\pi}{(1-p^2)^2}\right]^{1-1/K} \left(\int_{E^{\prime}}
w_0(z)^K J_{g_0}(z)\,dm \right)^{1/K}.
\ee
Now putting $w(z)=w_0(1/z),f(z)=g(1/z)$ and observing that
$J_g(z)=J_f(1/z)|z|^{-4},$ $J_{g_0}(z)=J_{f_0}(1/z)|z|^{-4},$ second inequality of
\eqref{eq13} follows from above. Here $E^{\prime}$ and $J_{g_0}$ is replaced by $E$
and $J_{f_0}=J_0$ respectively, where  $f_0(z)=1/(z-p).$  To obtain the first inequality
we use the other part of the `Harnack's Inequality' in \eqref{eq16} viz.
$$
h_p(\lambda) \geq  (1+|\lambda|)(1-|\lambda|)^{-1}h_p(0)
$$
and proceed in a similar fashion.
Next we show that the  second inequality of Theorem \ref{th2} is sharp.
To verify this, it is sufficient to show that the inequality \eqref{eq16a} is sharp. We follow the arguments given in \cite[Example 2.1]{AN}. First we
choose the numbers $w_j,p_j,r_j,\rho_j$ for $j=1,...,n$, suitably as
$1\leq w_1<w_2<\cdots<w_n$ and $0<p_j<1$, such that
\be\label{eq16b}
w_j=\left(\prod\limits_{l=1}^{j}r_l \right)^{-2/K} \quad\text{and}
\quad \sum\limits_{j=1}^{n} p_j w_j^K = 1.
\ee
We now consider the function
\be\label{eq17}
g=f^{\rho_1}_{r_1} \circ \cdots \circ f^{\rho_n}_{r_n},~~\text{where}~~
f^{\rho_j}_{r_j}(z)=\rho_jf_{r_j}(z/ \rho_j),\,j=1,...,n,
\ee
and $f_r$ defined in \eqref{eq01}. Next we consider the weight function
$w_0(z)= \sum\limits_{j=1}^{n} w_j\chi_{E_j}(z)$, where
$$
E_j=\lbrace z: \rho_{j+1}<|z|<\rho_jr_j \rbrace, ~~ 1\leq j\leq n-1; \quad E_n=
\lbrace z:|z|<\rho_nr_n\rbrace.
$$
The composition in \eqref{eq17} is well defined as we have
$$
r_j^2\rho_j^2-\rho_{j+1}^2=p_j,\,1\leq j\leq n-1;\,\, r_n^2\rho_n^2=p_n.
$$
In our case, we define
\be\label{eq18}
G(z)=(1-p^2)^{-1}g \left(\frac{z-p}{1-pz}\right) + p/(1-p^2),~~z\in \IC,
\ee
and the weight function as
$$
W_0(z)= \sum\limits_{j=1}^{n} w_j \chi_{\tilde{E_j}}(z), ~~ \tilde{E_j} =
\tilde{f}^{-1}(E_j),~~\text{where}~~
\tilde{f}(z)=(z-p)/(1-pz).
$$
Now the function $G$ defined in \eqref{eq18} belongs to the class
$\Sigma_k^0(p)$, as the function $g$ defined in \eqref{eq17} belongs to the
class $\Sigma_k^0$. If we now take $\tilde{E}=\cup _{j=1}^{n} \tilde{E_j}$,
then $G$ is conformal on $\tilde{E}$. Hence using first relation of \eqref{eq16b},
it is easy to see that
$$
J_{G|_{\tilde{E}}}= W_0(z)^{K-1}|1-pz|^{-4}= W_0(z)^{K-1}|g_0(\tilde{E})|,~z\in \tilde{E}.
$$
Again, using second relation of \eqref{eq16b}, we get
\beqq
\int_{\tilde{E}} W_0(z)J_G(z)\,dm
&=& \sum\limits_{j=1}^{n} \left( w_j^K \int_{\tilde{E_j}} |1-pz|^{-4} \,dm \right) \\
&=& \sum\limits_{j=1}^{n} \left( w_j^K |g_0(\tilde{E_j})| \right)\\
&=& \pi(1-p^2)^{-2} \left[ \sum\limits_{j=1}^{n-1} w_j^K(r_j^2\rho_j^2-\rho_{j+1}^2)
+  w_n^K r_n^2 \rho_n^2 \right]\\
&=& \pi(1-p^2)^{-2} \sum\limits_{j=1}^{n} p_jw_j^K \\
&=& \pi(1-p^2)^{-2}
= \left[\pi(1-p^2)^{-2}\right]^{1-1/K} \left(\int_{\tilde{E}}
W_0(z)^K J_{g_0}(z)\,dm\right)^{1/K}.
\eeqq
As equality holds in \eqref{eq16a}, hence it also holds for the second inequality in \eqref{eq13}. Optimality of the other inequality in \eqref{eq13} can be established by
similar construction.
\epf

\br
(i)~ While proving Theorem C in \cite{AN}, the authors first
assumed $E$ to be an open set and then proved the theorem for a general set
$E\subset \ID$ by limiting sense. Same argument also can be applied to the proof of
Theorem \ref{th2}, but we omit the details.

\noindent (ii)~ If $w(z)=1$ for all $z\in E,$ then the second inequality of Theorem~\ref{th2} implies Theorem~\ref{th1}(i).

\noindent (iii)~ In Theorem C, we assumed $f \in \Sigma_k^0$ of the form \eqref{eq00} in $\IDs,$ as taken in \cite{AN}. But if we take
$f \in \Sigma_k^0$ of the form \eqref{eq000} (with $b_0=0$) in $\ID$
and $f$ is conformal on $E\subset \IDs,$ then Theorem C can be restated
as
$$
\pi^{1-K}\left(\int_{E} w(z)^{1/K}|z|^{-4} \,dm \right)^K
\leq \int_{E} w(z)J_f(z)\,dm \leq \pi^{1-1/K} \left(\int_{E}
w(z)^K|z|^{-4} \,dm \right)^{1/K}.
$$
This result coincides with Theorem \ref{th2} for $p=0$.
\er

As an application of Theorem \ref{th1}, we prove the next result. It deals
with the bounds of the Hilbert transform of the characteristic function of a set $E\subset \ID.$
\bthm\label{th3}
If $E\subset \ID,$ then
\be\label{eq18a}
\int_{\ID \setminus E} \frac{1}{|1-pz|^{2}} \left|H \left[\frac{\chi_E}
{(1-p\overline{z})^{2}}\right] \right|\,dm \leq |g_0(E)| \log \left(
\frac{\pi(1-p^2)^{-2}}{|g_0(E)|} \right),
\ee
where $g_0(z)=z/(1-pz),z\in \IC$. The inequality is sharp.
\ethm

\bpf
For any function $\mu$ with $|\mu|=1,$ supported in $\ID \setminus E,$
we define $\mu_{\lambda}(z)=\lambda \mu(z)$ for $\lambda \in \ID$ and
consider the corresponding family of quasiconformal mappings $g_{\lambda}$
in $\sphere$, with dilatation $\mu_{\lambda}.$ We also assume that the
functions $g_{\lambda}$ are normalized such that they belong to
the class $\Sigma^0(p)$, when restricted on $\IDs$, therefore each
function $g_{\lambda}$ belongs to the class $\Sigma^0_{|\lambda|}(p)$,
for each $\lambda \in \ID$. Now by the assumption each $g_{\lambda}$ is conformal
on $E$, which gives from \eqref{eq9a} that
\beq\label{eq19}
|g_{\lambda}(E)| &=& \int_{E}|\dz g_{\lambda}(z)|^2\,dm \nonumber\\
&=& |g_0(E)| + 2 \mathrm{Re} \int_{E}(1-p\overline{z})^{-2}
H[\dzb g_{\lambda}]\,dm + \int_{E}|H[\dzb g_{\lambda}]|^2\,dm .
\eeq
Now from \eqref{eq9}, $w=\dzb g_{\lambda}$ can be written as
\begin{equation}\label{eq20}
\dzb g_{\lambda} = \lambda \mu(1-pz)^{-2} + h_{\lambda}(z),
\end{equation}
where
$\|h_{\lambda}\|_2 \leq C|\lambda|^2$, $C$ is a constant.
Using above identity it is easy to see that
\begin{equation*}
\int_{E}|H[\dzb g_{\lambda}]|^2\,dm = O(|\lambda|^2)\, \mbox{as}\, \lambda \to 0.
\end{equation*}
Again from \eqref{eq20} we get
\beqq
&&  \mathrm{Re} \int_{E}(1-p\overline{z})^{-2}H[\dzb g_{\lambda}(z)]\, dm \\
&=&  \mathrm{Re} \int_{E} \lambda(1-p\overline{z})^{-2}H[\mu(1-pz)^{-2}] \,dm
+ O(|\lambda|^2)\,, \quad \lambda \to 0.
\eeqq
Now upon using the last two estimates obtained above, we get from \eqref{eq19} that
\be\label{eq23}
|g_{\lambda}(E)|=|g_0(E)| + 2 \mathrm{Re} \int_{E} \lambda(1-p\overline{z})^{-2}
H[\mu (1-pz)^{-2}]\,dm + O(|\lambda|^2).
\ee

Now as $g_{\lambda} \in \Sigma^0_{|\lambda|}(p)$, by area distortion inequality
(Theorem \ref{th1}(i)), we get
$$
|g_{\lambda}(E)| \leq  \left[\pi(1-p^2)^{-2}\right]^{1-1/K} |g_0(E)|^{1/K},
$$
where $K=(1+|\lambda|)(1-|\lambda|)\inv$. Since
$1-K\inv=2|\lambda| + O(|\lambda|^2)$, therefore the above inequality can be written as
\beqq
|g_{\lambda}(E)| \leq
 |g_0(E)|+2|\lambda| |g_0(E)|\log \left(\pi(1-p^2)^{-2}|g_0(E)|^{-1} \right)
 + O(|\lambda|^2).
\eeqq
Comparing the coefficients of the terms which are linear in $|\lambda|$
of the above inequality and that of with \eqref{eq23}, we get
\begin{equation*}
 \mathrm{Re} \left[ \lambda \int_{E} (1-p\overline{z})^{-2} H[\mu (1-pz)^{-2}]\,dm \right] \leq |\lambda||g_0(E)| \log \left(\pi(1-p^2)^{-2}|g_0(E)|^{-1} \right).
\end{equation*}
Now for a particular choice of $\lambda $, we have
\begin{equation*}
 \mathrm{Re} \left[ \lambda \int_{E} (1-p\overline{z})^{-2} H[\mu (1-pz)^{-2}]\,dm \right] = |\lambda|\left| \int_{E} (1-p\overline{z})^{-2} H[\mu (1-pz)^{-2}]\,dm \right|.
\end{equation*}
From above two relations, we get
\be\label{eq24}
\left| \int_{E} (1-p\overline{z})^{-2} H[\mu (1-pz)^{-2}]\,dm \right|
\leq |g_0(E)| \log \left(\pi(1-p^2)^{-2}|g_0(E)|^{-1} \right).
\ee
Next using `symmetric property' of `$H$' (see \cite[p.95]{AIM}), we have
\beqq
\int_{E} (1-p\overline{z})^{-2} H[\mu (1-pz)^{-2}]\,dm
&=&  \int_{\IC} \chi_E (1-p\overline{z})^{-2} H[\mu (1-pz)^{-2}] \,dm \\
&=&  \int_{\IC} \mu(1-pz)^{-2} H[ \chi_E (1-p\overline{z})^{-2}]\,dm\\
&=&  \int_{\ID \setminus E} \mu(1-pz)^{-2} H[ \chi_E (1-p\overline{z})^{-2}]\,dm,
\eeqq
since $\mu$ has support in $\ID \setminus E.$ Using the inequality \eqref{eq24}, we get
\begin{equation*}
\left|\int_{\ID \setminus E} \mu(1-pz)^{-2} H[ \chi_E (1-p\overline{z})^{-2}]\,dm\right|
\leq |g_0(E)| \log \left(\pi(1-p^2)^{-2}|g_0(E)|^{-1} \right).
\end{equation*}
For a suitable choice of $\mu$, we can take modulus inside the integral of the
left hand side of the above inequality, which proves the theorem.
Finally it remains to prove the sharpness of the inequality \eqref{eq18a}. To show this we consider
$$
E= \left\lbrace z : \left|z-\frac{p(1-r^2)}{1-p^2r^2}\right|<\frac{r(1-p^2)}
{1-p^2r^2}\right\rbrace, \quad 0<r<1.
$$
Clearly $E\subset \ID$. Hence $|g_0(E)| = \pi r^2(1-p^2)^{-2},$ so that right hand side of \eqref{eq18a}
reduces to  $2\pi (1-p^2)^{-2}r^2 \log(r\inv)$. Next in order to find the Hilbert transform
of the function $\chi_E(1-p\overline{z})^{-2}$, we define
$$
f(z)=
\begin{cases}
\frac{1}{1-p^2}\left(\frac{\overline{z}-p}{1-p\overline{z}}\right), & z\in E,\\
\frac{r^2}{1-p^2}\left(\frac{1-pz}{z-p} \right), & z\in \IC \setminus E.
\end{cases}
$$
Here $f$ is continuous on $\IC$ and a little calculation reveals that $\dzb f=\chi_E(1-p\overline{z})^{-2}$ and
$\dz f = -r^2(z-p)^{-2}\chi_{\IC\setminus E}$. Using the relation
$H[\dzb f]=\dz f$, we have
$$
H[\chi_E(1-p\overline{z})^{-2}] = -r^2(z-p)^{-2}\chi_{\IC\setminus E}.
$$
Let $w=\tilde{f}(z)=(z-p)/(1-pz)=u+iv.$ Therefore, $\tilde{f}(\ID \setminus E)=
\lbrace w: r \leq |w| < 1 \rbrace $ and $J_{\tilde{f}}(z)=(1-p^2)^2|1-pz|^{-4}.$
Hence we have,
\beqq
\int_{\ID \setminus E} |1-pz|^{-2} \left|H[\chi_E(1-p\overline{z})^{-2}] \right|\,dm&=& r^2\int_{\ID\setminus E}
(|1-pz||z-p|)^{-2}\,dm\\
&=&  r^2(1-p^2)^{-2} \int_{\ID\setminus E} \left|\frac{1-pz}{z-p}\right|^2
\frac{(1-p^2)^2}{|1-pz|^4}\,dm\\
&=& r^2(1-p^2)^{-2} \int_{\tilde{f}(\ID\setminus E)}|w|^{-2} \,dudv \\
&=& 2\pi (1-p^2)^{-2}r^2 \log(r\inv).
\eeqq
Thus the inequality \eqref{eq18a} is sharp and this completes proof of the theorem.
\epf

\br
For $p=0$, the functions $g_{\lambda}$ defined in the proof of Theorem \ref{th3}
belong to the class $\Sigma^0_{|\lambda|}$ and the function $g_0$ becomes the identity function. Hence the inequality \eqref{eq18a} reads as (compare Theorem
14.6.1 of \cite[p.385]{AIM})
$$
\int_{\ID\setminus E} |H[\chi_E]| \,dm \leq |E|\log \left(\pi/|E|\right).
$$
\er

\end{document}